\documentstyle[12pt]{article}

\begin{document}

\begin{titlepage}

\title{Higher analytic torsion and cohomology of diffeomorphism groups}

\author{Ulrich Bunke\thanks{Mathematisches Institut, Universit\"at G\"ottingen,
Bunsenstr. 3-5, 37073 G\"ottingen, Germany, bunke@uni-math.gwdg.de}}
\end{titlepage}

\newcommand{\diag}{{\rm diag}}
\newcommand{\proof}{{\it Proof.$\:\:\:\:$}}
\newcommand{\Bbb}{\rm}
\newcommand{\dist}{{\rm dist}}
\newcommand{\kaaa}{{\bf k}}
\newcommand{\paaa}{{\bf p}}
\newcommand{\taaa}{{\bf t}}
\newcommand{\haaa}{{\bf h}}
\newcommand{\R}{{\bf R}}
\newcommand{\Q}{{\bf Q}}
\newcommand{\Z}{{\bf Z}}
\newcommand{\C}{{\bf C}}
\newcommand{\K}{{\tt K}}
\newcommand{\Naaa}{{\bf N}}
\newcommand{\gaaa}{{\bf g}}
\newcommand{\maaa}{{\bf m}}
\newcommand{\aaaa}{{\bf a}}
\newcommand{\naaa}{{\bf n}}
\newcommand{\brr}{{\bf r}}
\newcommand{\res}{{\rm res}}
\newcommand{\Tr}{{\rm Tr}}
\newcommand{\cT}{{\cal T}}
\newcommand{\dom}{{\rm dom}}
\newcommand{\db}{{\bar{\partial}}}
\newcommand{\g}{{\gaaa}}
\newcommand{\cZ}{{\cal Z}}
\newcommand{\cH}{{\cal H}}
\newcommand{\cM}{{\cal M}}
\newcommand{\interi}{{\rm int}}
\newcommand{\singsupp}{{\rm singsupp}}
\newcommand{\cE}{{\cal E}}
\newcommand{\cV}{{\cal V}}
\newcommand{\cI}{{\cal I}}
\newcommand{\cC}{{\cal C}}
\newcommand{\mod}{{\rm mod}}
\newcommand{\cK}{{\cal K}}
\newcommand{\cA}{{\cal A}}
\newcommand{\cEp}{{{\cal E}^\prime}}
\newcommand{\cU}{{\cal U}}
\newcommand{\Hom}{{\mbox{\rm Hom}}}
\newcommand{\vol}{{\rm vol}}
\newcommand{\cO}{{\cal O}}
\newcommand{\End}{{\mbox{\rm End}}}
\newcommand{\Ext}{{\mbox{\rm Ext}}}
\newcommand{\rk}{{\mbox{\rm rank}}}
\newcommand{\im}{{\mbox{\rm im}}}
\newcommand{\sign}{{\rm sign}}
\newcommand{\spann}{{\mbox{\rm span}}}
\newcommand{\symm}{{\mbox{\rm symm}}}
\newcommand{\cF}{{\cal F}}
\newcommand{\Ree}{{\rm Re }}
\newcommand{\Res}{{\mbox{\rm Res}}}
\newcommand{\Imm}{{\rm Im}}
\newcommand{\inter}{{\rm int}}
\newcommand{\clo}{{\rm clo}}
\newcommand{\tg}{{\rm tg}}
\newcommand{\ee}{{\rm e}}
\newcommand{\Li}{{\rm Li}}
\newcommand{\cN}{{\cal N}}
 \newcommand{\conv}{{\rm conv}}
\newcommand{\op}{{\mbox{\rm Op}}}
\newcommand{\tr}{{\mbox{\rm tr}}}
\newcommand{\cs}{{c_\sigma}}
\newcommand{\ctg}{{\rm ctg}}
\newcommand{\degg}{{\mbox{\rm deg}}}
\newcommand{\Ad}{{\mbox{\rm Ad}}}
\newcommand{\ad}{{\mbox{\rm ad}}}
\newcommand{\codim}{{\mbox{\rm codim}}}
\newcommand{\Gr}{{\mbox{\rm Gr}}}
\newcommand{\coker}{{\rm coker}}
\newcommand{\id}{{\mbox{\rm id}}}
\newcommand{\ord}{{\mbox{\rm ord}}}
\newcommand{\nat}{{\bf  N}}
\newcommand{\supp}{{\mbox{\rm supp}}}
\newcommand{\spec}{{\mbox{\rm spec}}}
\newcommand{\Ann}{{\mbox{\rm Ann}}}
\newcommand{\aca}{{\aaaa_\C^\ast}}
\newcommand{\acag}{{\aaaa_{\C,good}^\ast}}
\newcommand{\acage}{{\aaaa_{\C,good}^{\ast,extended}}}
\newcommand{\ck}{{\cal K}}
\newcommand{\tck}{{\tilde{\ck}}}
\newcommand{\tnk}{{\tilde{\ck}_0}}
\newcommand{\ceep}{{{\cal E}(E)^\prime}}
 \newcommand{\ncE}{{{}^\naaa\cE}}
 \newcommand{\Or}{{\rm Or }}

\newcommand{\cB}{{\cal B}}
\newcommand{\hc}{{{\cal HC}(\gaaa,K)}}
\newcommand{\hcma}{{{\cal HC}(\maaa_P\oplus\aaaa_P,K_P)}}

\newcommand{\vsl}{{V_{\sigma_\lambda}}}
\newcommand{\czg}{{\cZ(\gaaa)}}
\newcommand{\csl}{{\chi_{\sigma,\lambda}}}
\newcommand{\cR}{{\cal R}}
\def\hB{\hspace*{\fill}$\Box$ \newline\noindent}
\newcommand{\varho}{\varrho}
\newcommand{\ind}{{\rm index}}
\newcommand{\Ind}{{\rm Ind}}
\newcommand{\Fin}{{\mbox{\rm Fin}}}
\newcommand{\cS}{{\cal S}}
\newcommand{\orig}{{\cal O}}
\def\hB{\hspace*{\fill}$\Box$ \\[0.5cm]\noindent}
\newcommand{\cL}{{\cal L}}
 \newcommand{\cG}{{\cal G}}
\newcommand{\npci}{{\naaa_P\hspace{-1.5mm}-\hspace{-2mm}\mbox{\rm coinv}}}
\newcommand{\pki}{{(\paaa,K_P)\hspace{-1.5mm}-\hspace{-2mm}\mbox{\rm inv}}}
\newcommand{\mki}{{(\maaa_P\oplus \aaaa_P, K_P)\hspace{-1.5mm}-\hspace{-2mm}\mbox{\rm inv}}}

\newcommand{\npi}{{\naaa_P\hspace{-1.5mm}-\hspace{-2mm}\mbox{\rm inv}}}
\newcommand{\ngp}{{N_\Gamma(\pi)}}
\newcommand{\gbg}{{\Gamma\backslash G}}
\newcommand{\gkm}{{ Mod(\gaaa,K) }}
\newcommand{\ggkm}{{  (\gaaa,K) }}
\newcommand{\pkm}{{ Mod(\paaa,K_P)}}
\newcommand{\ppkm}{{  (\paaa,K_P)}}
\newcommand{\makm}{{Mod(\maaa_P\oplus\aaaa_P,K_P)}}
\newcommand{\mmakm}{{ (\maaa_P\oplus\aaaa_P,K_P)}}
\newcommand{\cP}{{\cal P}}
\newcommand{\gm}{{Mod(G)}}
\newcommand{\gk}{{\Gamma_K}}
\newcommand{\La}{{\cal L}}
\newcommand{\cug}{{\cU(\gaaa)}}
\newcommand{\cuk}{{\cU(\kaaa)}}
\newcommand{\dc}{{C^{-\infty}_c(G) }}
\newcommand{\gdk}{{\gaaa/\kaaa}}
\newcommand{\dgkm}{{ D^+(\gaaa,K)-\mbox{\rm mod}}}
\newcommand{\dgm}{{D^+G-\mbox{\rm mod}}}
\newcommand{\vect}{{\C-\mbox{\rm vect}}}
 \newcommand{\cig}{{C^{ \infty}(G)_{K} }}
\newcommand{\gami}{{\Gamma\hspace{-1.5mm}-\hspace{-2mm}\mbox{\rm inv}}}
\newcommand{\cQ}{{\cal Q}}
\newcommand{\mmap}{{Mod(M_PA_P)}}

\newtheorem{prop}{Proposition}[section]
\newtheorem{lem}[prop]{Lemma}
\newtheorem{ddd}[prop]{Definition}
\newtheorem{theorem}[prop]{Theorem}
\newtheorem{kor}[prop]{Corollary}
\newtheorem{ass}[prop]{Assumption}
\newtheorem{con}[prop]{Conjecture}
\newtheorem{prob}[prop]{Problem}
\newtheorem{fact}[prop]{Fact}

\maketitle

\begin{abstract}
We consider a closed odd-dimensional oriented manifold $M$ together with an acyclic
flat hermitean vector bundle $\cF$. We form the trivial fibre bundle
with fibre $M$ over the manifold of all Riemannian metrics
on $M$. It has a natural flat connection and a vertical Riemannian
metric. The higher analytic torsion form of Bismut/Lott
associated to the situation is invariant with respect to
the connected component of the identity of the diffeomorphism group of $M$.
Using that the space of Riemannian metrics is contractible
we define continuous cohomology classes of the diffeomorphism group
and its Lie algebra.
For the circle we compute this classes in degree 2
and show that the group cohomology class is non-trivial, while
the Lie algebra cohomology class vanishes.
\end{abstract}

\tableofcontents

\newcommand{\Met}{{{\cal M}et}}
\newcommand{\Diff}{{{\cal D}iff}}
\newcommand{\cX}{{\cal X}}
\newcommand{\cD}{{\cal D}}

\section{The form T}\label{tt}

Let $M$ be a closed odd-dimensional oriented manifold.
By $\Met(M)$ we denote the space of all Riemannian metrics
on $M$. It is a convex open subset of $C^\infty(M,S^2T^*M)$
and therefore a Fr\'echet manifold.

By $\Diff(M)$ we denote the group of
diffeomorphisms of $M$. It is a Fr\'echet Lie group,
and its Lie algrabra is $\cX(M)=C^\infty(M,TM)$.
The group $\Diff(M)$ acts on $\Met(M)$ by
$$(f,g)\in \Diff(M)\times \Met(M)\mapsto w(f)g:=(f^{-1})^* g\in \Met(M)\ .$$
By $\Diff(M)^0$ we denote the connected component
of the identity of $\Diff(M)$.

Let $F\rightarrow M$ be a flat hermitan vector bundle.
By $\Diff(M,F)$ we denote the group of its automorphisms.
$\Diff(M,F)$ is again a Fr\'echet Lie group, and
we denote by $\Diff(M,F)^0$ its connected component of the identity.
There is a natural surjection
$q_F:\Diff(M,F)^0\rightarrow \Diff(M)^0$.

\begin{ass}\label{acy}
We assume that $\cF$ is acyclic, i.e., that $H^*(M,\cF)=0$,
where $\cF$ is the sheaf of parallel sections of $F$.
\end{ass}

We are going to define a closed and $\Diff(M)^0$-invariant form
$$T=T_0+T_2+T_4\dots, \quad T_{2i}\in \Omega^{2i}(\Met(M))$$
by specializing the higher real analytic torsion
introduced by Bismut/Lott \cite{bismutlott95}.

Let $\Omega^*(M,F)$ denote the space of $F$-valued forms on $M$.
If we choose a Riemannian metric $g\in \Met(M)$,
then we can define an $L^2$-scalar product on $\Omega^*(M,F)$.
For $X\in\cX(M)$ let $i(X)\in\End(\Omega^*(M,F))$ denote
the insertion of $X$. We put $e(X):=i(X)^*$ and define
$$c(X):=e(X)-i(X), \quad \hat{c}(X):=e(X)+i(X)\ .$$
Let $\nabla$ denote the connection on $\Lambda^*T^*M\otimes F$
induced by the Levi-Civita connection associated to
$g$ and the flat  connection on $F$.

We identify $T_g\Met(M)\cong C^\infty(M,S^2T^*M)$ in the natural way.
Then we define $S\in\Hom(T_g\Met(M),\End(\Omega^*(M,F)))$
by $$S(h):=h(e_i,e_j) c(e_i)\hat{c}(e_j),\quad h\in T_g\Met(M)\ ,$$
where $\{e_i\}$ denotes a local orthonormal frame on $M$.

We put $\cA^k:= \Hom(\Lambda^kT_g\Met(M),\End(\Omega^*(M,F)))$ and
$\cA:=\oplus_k \cA^k$.
Then there is a natural product
$\cA^k\otimes \cA^l\rightarrow \cA^{k+l}$.

For $t > 0$ we define
\begin{equation}\label{yuu}D_t:=-\sqrt{t}\hat{c}(e_i)\nabla_{e_i}+S\in \cA\ .\end{equation}
By $N\in \cA^0$ we denote the $\Z$-grading on $\Omega^*(M,F)$.
We have $D_t^2:=-t\Delta \mbox{(mod $\cA^{>0}$)}$,
where $\Delta$ is the Laplacian on $\Omega^*(M,F)$ defined with the Riemannian
metric $g$.
Hence for any $i\in\nat$ and for $t > 0$ the term 
$[(1+2D_t^2)\ee^{D_t^2}]_{i}\in \cA^i$
has values in the trace class operators on $\Omega(M,F)$.

\begin{lem}
For all $i\in 2\nat_0$ and $h\in \Lambda^{i}T_g\Met(M)$
the integral
\begin{equation}\label{lll1} T_i(h):=-(\frac{1}{2\pi\imath})^{i/2}\int_0^\infty \Tr_s N [(1+2D^2 _t) \ee^{D_t^2}]_i(h)
                \frac{dt}{t} \end{equation}
converges. Moreover
$T:=T_0+T_2+T_4+\dots$ is a closed form in $\Omega^{ev}(\Met(M))$.
\end{lem}
\proof
We show how this can be deduced from the results of \cite{bismutlott95}
by specialization. We consider the trivial fibre bundle
$p:E:=\Met(M)\times M\rightarrow \Met(M)$. Let $q:E\rightarrow M$
denote the projection onto the second factor.
Then $q^* F\rightarrow E$
is a flat hermitean vector bundle over $E$. We choose the tautological
vertical Riemannian metric on $E$ such that the fibre $E_g$ carries
the metric $g\in\Met(M)$. There is a natural choice of a horizontal
distribution $T^H E\rightarrow E$ given by the kernel of $dq$.

The tensor $T$ given in \cite{bismutlott95}, (3.11), vanishes since $T^H E$ is
integrable. Since $q^*F$ is flat as a hermitean vector bundle
the tensor $\psi$ introduced in \cite{bismutlott95}, (d),
vanishes, too. The symmetric tensor $\omega_{\alpha k j}$ defined
in \cite{bismutlott95}, (3.21), is can be identified with the linear
map $h\in T_g\Met(M)\rightarrow \omega(h)_{i j}=h(e_i,e_j)$.
Thus $D_1$ coincides with $2X$ defined in \cite{bismutlott95}, (3.41).
One can now check that $D_t$ is just  twice the operator
given in \cite{bismutlott95}, (3.50).

Convergence of the integral (\ref{lll1}) follows from \cite{bismutlott95}, Thm. 3.21.
By \cite{bismutlott95}, Cor. 3.25  the form $T$ is closed. \hB

\begin{lem}
The form $T$ is invariant under $\Diff(M)^0$.
\end{lem}
\proof
Via the homomorphism $q_F:\Diff(M,F)^0\rightarrow  \Diff(M)^0$
the  group $\Diff(M,F)^0$
acts on $M$, $\Met(M)$, and thus on $E$
such that the projections $p,q$ are equivariant.
Since $q$ is equivariant, the horizontal distribution
$T^HE$ is $\Diff(M,F)^0$ invariant.
The tautological vertical Riemannian metric
is invariant with respect to $\Diff(M,F)^0$, too.
Now $\Diff(M,F)^0$  acts on $q^*F\rightarrow E$ by automorphisms
of flat hermitean vector bundles.

Since all structures used to
define $T$ are $\Diff(M,F)^0$-invariant, we conclude that $T$
is invariant with respect to $\Diff(M,F)^0$, too
(compare \cite{bismutlott95}. Thm. A 1.1).
Since $q_F:\Diff(M,F)^0\rightarrow \Diff(M)^0$ is surjective,
$T$ is invariant under $\Diff(M)^0$.
\hB

\section{Group cohomology classes of $\Diff(M)^0$}

Let $T$ be a closed $p$-form on the convex subset $\Met(M)$
of $C^\infty(M,S^2 T^*M)$ which is invariant under $\Diff(M)^0$.
We fix a base point $g_b\in\Met(M)$.
If $f_0,\dots f_p\in \Diff(M)^0$, then we define a smooth map
$$s(f_0,\dots,f_p):\Delta^p\rightarrow \Met(M)$$
from the standard $p$-simplex $\Delta^p$ to $\Met(M)$
by $s(f_0,\dots,f_p)(t):=\sum_{i=0}^p t_i w(f)g_b$,
where $t=(t_0,\dots,t_p)\in \Delta^p$, $\sum_{i=0}^p t_i = 1$.
If $f\in\Diff(M)^0$, then
$s(f f_0,\dots,f f_p)=w(f)s(f_0,\dots,f_p)$.

Let
$C^p(\Diff(M)^0)$ be
the space of real alternating group $p$-cochains,
and let $d:C^p(\Diff(M)^0)\rightarrow C^{p+1}(\Diff(M)^0)$
be the  usual differential \cite{brown82}, Ch 1.5.
Then $\Diff(M)^0$ acts on $C^p(\Diff(M)^0)$ by $(fc)(f_0,\dots,f_p):=c(f^{-1}f_0,\dots,f^{-1}f_p)$.
We define the cochain $c_T$ by
$$c_T(f_0,\dots,f_p):=\int_{\Delta^p} s(f_0,\dots,f_p)^* T\ .$$
\begin{lem}
$c_T$ is a $\Diff(M)^0$-invariant cocycle.
\end{lem}
\proof
Let $f\in Diff(M)^0$.
Then
\begin{eqnarray*}
(fc_T)(f_0,\dots,f_p)&=&c_T(f^{-1}f_0,\dots,f^{-1}f_p)\\
&=&\int_{\Delta^p} s(f^{-1}f_0,\dots,f^{-1}f_p)^* T\\
&=&\int_{\Delta^p} (w(f^{-1})\circ s(f_0,\dots,f_p))^* T\\
&=&\int_{\Delta^p} s(f_0,\dots,f_p)^* w(f^{-1})^*T\\
&=&\int_{\Delta^p} s(f_0,\dots,f_p)^* T\\
&=&c_T(f_0,\dots,f_p)\ .
\end{eqnarray*}
We have by Stokes Lemma
\begin{eqnarray*}
(dc_T)(f_0,\dots,f_{p+1})&=&\sum_{i=0}^{p+1}(-1)^p c_T(f_0,\dots,\hat{f}_i,\dots,f_{p+1})\\
&=&\sum_{i=0}^{p+1}(-1)^p \int_{\Delta^p} s(f_0,\dots,\hat{f}_i,\dots,f_{p+1})^* T\\
&=&\int_{\partial \Delta^{p+1}} s(f_0,\dots,f_{p+1})^* T\\
&=&\int_{\Delta^{p+1}} d (s(f_0,\dots,f_{p+1})^* T)\\
&=&\int_{\Delta^{p+1}} s(f_0,\dots,f_{p+1})^* dT\\
&=&0 \ .
\end{eqnarray*}\hB

Thus $c_T$ defines a cohomology class $h_T\in H^*(\Diff(M)^0,\R)$.
\begin{lem}
$h_T$ does not depend on the choice of the base point $g_b$.
\end{lem}
\proof
Let $g_b^\prime\in\Met(M)$ be another base point and $c_T^\prime$
be defined using $g_b^\prime$.
Then we define the chain $u\in C^{p-1}(\Diff(M)^0)$ by
$$u(f_0,\dots,f_{p-1})=\int_{\Delta^{p-1}\times I}
         U_{p-1}(f_0,\dots,f_{p-1})^* T\ , $$
where $U_{p-1}(f_0,\dots,f_{p-1}):\Delta^{p-1}\times I\rightarrow \Met(M)$
is given by
$U_{p-1}(f_0,\dots,f_{p-1})(t,s):=\sum_{i=0}^{p-1} t_i (s g_b+(1-s)g_b^\prime)$.
Now we have
\begin{eqnarray*}
c_T^\prime(f_0,\dots,f_p)-c_T(f_0,\dots,f_p) -du(f_0,\dots,f_p)
&=&\int_{\partial (\Delta^p\times I)} U_{p}(f_0,\dots,f_p)^*T\\
&=&\int_{\Delta^p\times I} U_{p}(f_0,\dots,f_p)^*dT\\
&=&0\ .
\end{eqnarray*}\hB
Recall that $\Diff(M)^0$ is a Fr\'echet Lie group.
Let $C_c^*(\Diff(M)^0)\subset C^*(\Diff(M)^0)$ denote the subcomplex of smooth
cochains and $H^*_c(\Diff(M),\R)$ be its cohomology.
Then  $u,c_T\in C_c^*(\Diff(M)^0)$, and $h_T\in  H^p_c(\Diff(M),\R)$
is well-defined independently of the choice of $g_b$.

\begin{ddd}In the special case that $T$ is the higher analytic torsion form 
associated to the closed odd-dimensional oriented manifold $M$ and the 
acyclic locally constant sheaf of Hilbert spaces $\cF$ defined
in Section \ref{tt} we denote the class $h_T$ by $\cT(M,\cF)$.
\end{ddd}

\section{Lie algebra cohomology classes of $\cX(M)$}

Recall that $\cX(M)$ is the Lie algebra of $Diff(M)^0$.
Let
$C^*(\cX(M)):=\Hom(\Lambda^{*}\cX(M),\R)$
denote the complex of continuous Lie algebra cochains with differential $d$
(see \cite{fuks84}, Ch. 1.3).
There is a natural map of cochain
complexes $\cD:C_c^*(\Diff(M)^0)\rightarrow C^*(\cX(M))$
given by
$$(\cD c)(X_1,\dots,X_p):=p\frac{d}{dt_1}_{t_1=0}\dots \frac{d}{dt_p}_{t_p=0}
c(1,\ee^{t_1X_1},\dots,\ee^{t_p X_p})\ ,$$
$c\in C^p_c(\Diff(M)^0)$, $X_i\in \cX(M)$.
Let $\cD_*:H_c(\Diff(M),\R)\rightarrow H^*(\cX(M),\R)$ denote the induced
map.

If $X\in \cX(M)$, then we define
$$h_X:=\frac{d}{dt}_{t=0} w(\ee^{tX})g_b\in T_{g_b}\Met(M)\ .$$
\begin{lem}\label{r33}
Let $T$ be a closed $p$-form on $\Met(M)$.
We have $$(\cD c_T)(X_1,\dots, X_p)=\frac{1}{(p-1)!}T(h_{X_1}\wedge\dots\wedge h_{X_p})\ .$$
\end{lem}
\proof
We have
\begin{eqnarray*}
s(1,\ee^{\epsilon_1 X_1},\dots,\ee^{\epsilon_p X_p})(t)&=&
t_0g_b + \sum_{i=1}^p t_i \ee^{\epsilon_i X_i}(g_b)\\
&=&t_0 g_b + \sum_{i=1}^p t_i (g_b+\epsilon_i h_{X_i}+O(\epsilon_i^2))\\
d s(1,\ee^{\epsilon_1 X_1},\dots,\ee^{\epsilon_p X_p})(t)&=
& g_b dt_0+\sum_{i=1}^p dt_i (g_b+\epsilon_i h_{X_i}+O(\epsilon_i^2))\\
&=&\sum_{i=0}^p dt_i   g_b + \sum_{i=1}^p dt_i
            (\epsilon_i h_{X_i}+O(\epsilon_i^2))\\
&=&\sum_{i=1}^p dt_i (\epsilon_i h_{X_i}+O(\epsilon_i^2))\\
s(1,\ee^{\epsilon_1 X_1},\dots,\ee^{\epsilon_p X_p})^* T(t)
&=&(\prod_{i=1}^p \epsilon_i )T(h_{X_1}\wedge\dots\wedge h_{X_p})
            (\prod_{i=1}^p dt_i) + O(\epsilon_1^2,\dots, \epsilon_p^2)\\
c_T(1,\ee^{\epsilon_1 X_1},\dots,\ee^{\epsilon_p X_p})&=&
(\prod_{i=1}^p \epsilon_i) T(h_{X_1}\wedge\dots\wedge h_{X_p})
   \int_{\Delta^p}(\prod_{i=1}^p dt_i) + O(\epsilon_1^2,\dots, \epsilon_p^2)\\
&=&\frac{1}{p!}T(h_{X_1}\wedge\dots\wedge h_{X_p})   (\prod_{i=1}^p \epsilon_i)
         + O(\epsilon_1^2,\dots, \epsilon_p^2)\ .
\end{eqnarray*}
This implies the Lemma. \hB

\section{Computations for the circle}
\newcommand{\Arg}{{\rm Arg}}

In this section we explicitly compute $T_2\in\Omega^2(\Met(S^1))$,
$\cT_2(M,\cF)\in H^2_c(\Diff(M),\R)$ and $\cD_*\cT_2(S^1,\cF) \in H^2(\cX(M),\R)$
in the case that $M=S^1$ and $\cF$ is a locally constant sheaf
of one-dimensional Hilbert spaces with non-trivial holonomy.

We identify $S^1\cong \R/\Z$.
Let $F$ be the flat hermitean vector bundle associated to $\cF$. 
Up to isomorphism $F$ is determined by its holonomy
$\ee^{2\pi\imath a}$, $a\in [0,1)$. If $a\in (0,1)$,
then $\cF$ is acyclic.

We make the identification $\Omega^p(S^1,F)\cong \{f\in C^\infty(\R)\:|\: f(x+1)=\ee^{2\pi\imath a}f(x)\}:=\cH$ for
$p=0,1$ using the form $dx$ in the case $p=1$.
We further identify $C^\infty (S^1,S^2T^*S^1)$ with $C^\infty(S^1)$
using the metric $dx^2$.

Let $g_b$ be the standard metric $dx^2$ on $S^1$ of volume $1$.
Then $w(\Diff(S^1)^0 )dx^2=\Met_1(S^1)$, where
$\Met_1 (S^1)=\{g\in \Met(S^1)\:|\: \vol_g(S^1)=1\}$.
In order to compute the restriction of $T$ to $\Met_1(S^1)$
it is therefore sufficient to compute $T(dx^2)$.

In the following we fix the metric $dx^2$.
We put $i:=i(\partial_x)$, $e:=e(\partial_x)=dx \wedge$. Then
$c=e-i$, $\hat{c}=e+i$, and $c\hat{c}=-z$, where 
$z$ is the $\Z_2$-grading of $\Omega^*(S^1,F)$.

Let $V\subset T_{dx^2}\Met(S^1)$ be a finite-dimensional
subspace and $\Gr(V^*)$ be the Grassmann algebra generated by
$V^*$. Then $\Gr(V^*)$ is a $\Z$-graded algebra,
and we denote by $\Gr(V^*)_p$ the subspace of elements of degree $p$.
We choose a base $\{h_\alpha\}$ of $V$ and let $\{E^\alpha\}$
denote the dual base.
Then we define
$$D:=-\hat{c}\partial_x - h_\alpha E^\alpha z$$
acting on $\cH\otimes \C^2$, where
$$\hat{c}:=\left(\begin{array}{cc}0&1\\1&0\end{array}\right),\quad
z=\left(\begin{array}{cc}1&0\\0&-1\end{array}\right)\ ,$$
$h_\alpha\in C^\infty(S^1)$ acts as multiplication operator  on $\cH$,
and $\{E^\alpha,\hat{c}\}=0=[E^\alpha,z]$.

The operator $D$ corresponds to $D_1$ in (\ref{yuu}) under the
identifications above. We obtain $D_t$ by rescaling.
For $t>0$ let $\Psi_t:V\rightarrow V$ be multiplication by $\frac{1}{\sqrt{t}}$
and $\Psi_t^*:\Gr(V^*)\rightarrow \Gr(V^*)$ be the induced automorphism.
Then we have $D_t=\Psi_t^* \sqrt{t}D$.

We have
\begin{eqnarray*}
D^2&=&\partial_x^2 + \{\hat{c}\partial_x, h_\alpha E^\alpha z\} + h_\alpha h_\beta E^\alpha E^\beta\\
&=&\partial_x^2 + 2 \hat{c} h_\alpha E^\alpha z\partial_x+ \hat{c}(\partial_x h_\alpha) E^\alpha z\ .
\end{eqnarray*}
We can write $\Tr_s N(1+2D^2_t)\ee^{D_t^2}=\Psi^*_t(1+2t\frac{d}{dt}) \Tr_s N\ee^{tD^2}$,
where 
$$N=\left(\begin{array}{cc}0&0\\0&1\end{array}\right) \ . $$
We want to compute $[\ee^{tD^2}]_2\in \Gr(V^*)_2\otimes \End(\Omega^*(S^1,F))$.
Let $R_\alpha:=  -2 \hat{c} h_\alpha  z\partial_x - \hat{c}(\partial_x h_\alpha)  z$.
Then we have $D^2= \partial_x^2+ E^\alpha R_\alpha $.
By Duhamel's formula we have
\begin{eqnarray*}
[\ee^{tD^2}]_2&=& -t^2 E^\alpha E^\beta \int_{\Delta^2}\ee^{t\sigma_0\partial_x^2} R_\alpha \ee^{t\sigma_1\partial_x^2} R_\beta \ee^{t\sigma_2\partial_x^2} d\sigma\\
&=&t^2 E^\alpha E^\beta \int_{\Delta^2}\ee^{t\sigma_0\partial_x^2} (2h_\alpha\partial_x+\partial_xh_\alpha) \ee^{t\sigma_1\partial_x^2}
(2h_\beta\partial_x+\partial_xh_\beta) \ee^{t\sigma_2\partial_x^2} d\sigma\ ,\\
\Tr_s N [\ee^{tD^2}]_2&=&- t^2 E^\alpha E^\beta \Tr \int_{\Delta^2}\ee^{t\sigma_0\partial_x^2} (2h_\alpha\partial_x+\partial_xh_\alpha) \ee^{t\sigma_1\partial_x^2}
(2h_\beta\partial_x+\partial_xh_\beta) \ee^{t\sigma_2\partial_x^2} d\sigma\ ,
\end{eqnarray*}
where the operator in the last line acts on $\cH$.

For $k\in\Z$ let  $f_k(x):=\ee^{2\pi\imath (k+a)x}$. Then $\partial_x^2 f_k=-4\pi^2 (k+a)^2 f_k$. 
Furthermore, for $\alpha\in\Z$ let $h_\alpha(x):=\ee^{2\pi\imath \alpha x}$.
Let $V$ be spanned by a finite number of these $h_\alpha$.
Then we have
\begin{eqnarray*}
\ee^{t\sigma_2\partial_x^2}f_k&=&\ee^{-4\pi^2 (k+a)^2 t \sigma_2} f_k\\
(2h_\beta\partial_x+\partial_xh_\beta) \ee^{t\sigma_2\partial_x^2}f_k&=&
2\pi\imath (2(k+a)+\beta)   \ee^{-4\pi^2 (k+a)^2 t \sigma_2}f_{k+\beta}\\&&\hspace{-3cm}
 \ee^{t\sigma_1\partial_x^2}
(2h_\beta\partial_x+\partial_xh_\beta) \ee^{t\sigma_2\partial_x^2}f_k\\&=&
2\pi\imath \ee^{-4\pi^2 (k+a+\beta)^2 t \sigma_1} (2(k+a)+\beta)   \ee^{-4\pi^2 (k+a)^2 t \sigma_2}f_{k+\beta}\\&&\hspace{-3cm}
(2h_\alpha\partial_x+\partial_xh_\alpha) \ee^{t\sigma_1\partial_x^2}
(2h_\beta\partial_x+\partial_xh_\beta) \ee^{t\sigma_2\partial_x^2}\\&=&
 -4\pi^2 (2(k+a+\beta)+\alpha)\\&& \ee^{-4\pi^2 (k+a+\beta)^2 t \sigma_1} (2(k+a)+\beta)    \ee^{-4\pi^2 (k+a)^2 t \sigma_2}f_{k+\beta+\alpha}\\&&\hspace{-3cm}
\ee^{t\sigma_0\partial_x^2} (2h_\alpha\partial_x+\partial_xh_\alpha) \ee^{t\sigma_1\partial_x^2}
(2h_\beta\partial_x+\partial_xh_\beta) \ee^{t\sigma_2\partial_x^2}f_k\\&=&
 -4\pi^2 \ee^{-4\pi^2 (k+a+\beta+\alpha)^2 t \sigma_0}(2(k+a+\beta)+\alpha)\\&& \ee^{-4\pi^2 (k+a+\beta)^2 t \sigma_1} (2(k+a)+\beta)   \ee^{-4\pi^2 (k+a)^2 t \sigma_2}f_{k+\beta+\alpha}
\end{eqnarray*}
We conclude that if $\alpha+\beta\not=0$, then
$$\Tr\:  \ee^{t\sigma_0\partial_x^2} (2h_\alpha\partial_x+\partial_xh_\alpha) \ee^{t\sigma_1\partial_x^2}
(2h_\beta\partial_x+\partial_xh_\beta) \ee^{t\sigma_2\partial_x^2} =0\ .$$
If $\alpha=-\beta$, then
\begin{eqnarray*}
&&\hspace{-2cm}\Tr \:\ee^{t\sigma_0\partial_x^2} (2h_\alpha\partial_x+\partial_xh_\alpha) \ee^{t\sigma_1\partial_x^2}
(2h_\beta\partial_x+\partial_xh_\beta) \ee^{t\sigma_2\partial_x^2} 
\\&=&-4\pi^2 \sum_{k\in\Z}
 \ee^{-4\pi^2 (k+a)^2 t \sigma_0}(2(k+a)-\alpha) \ee^{-4\pi^2 (k+a-\alpha)^2 t \sigma_1} (2(k+a)-\alpha)    \ee^{-4\pi^2 (k+a)^2 t\sigma_2}\\
&=&-4\pi^2 \sum_{k\in\Z}(2(k+a)-\alpha)^2\ee^{-4\pi^2 (k+a)^2 t}
\ee^{4\pi^2\alpha(2(k+a)-\alpha) t \sigma_1}\ .
\end{eqnarray*}
Now 
\begin{eqnarray*}
\int_{\Delta^2} f(\sigma_1) d\sigma &=&
 \int_0^1 f(\sigma_1) \int_0^{\sigma_1} d\sigma_2d\sigma_1  \\
&=&\int_0^1 f(\sigma_1)\sigma_1 d\sigma_1\ .
\end{eqnarray*}
We conclude 

\begin{eqnarray*}
\lefteqn{\Tr\:  \int_{\Delta^2}  \ee^{t\sigma_0\partial_x^2} (2h_\alpha\partial_x+\partial_xh_\alpha) \ee^{t\sigma_1\partial_x^2}
(2h_\beta\partial_x+\partial_xh_\beta) \ee^{t\sigma_2\partial_x^2} d\sigma}\\
 &=& -4\pi^2 \sum_{k\in\Z} \int_0^1 
 (2(k+a)-\alpha)^2\ee^{-4\pi^2 (k+a)^2 t}
\ee^{4\pi^2\alpha(2(k+a)-\alpha) t \sigma}
\sigma d\sigma\\
&=& -4\pi^2 \sum_{k\in\Z}\ee^{-4\pi^2 (k+a)^2 t} \int_0^{2(k+a)-\alpha}\ee^{4\pi^2\alpha t u} u du\\
&=& -4\pi^2 \sum_{k\in\Z} \ee^{-4\pi^2 (k+a)^2 t}
\left(\frac{\ee^{4\pi^2 \alpha t(2(k+a)-\alpha)}(2(k+a)-\alpha)}{4\pi^2\alpha t}-\frac{\ee^{4\pi^2 \alpha t(2(k+a)-\alpha)}-1}{(4\pi^2\alpha t)^2}\right)\\
&=& -4\pi^2 \sum_{k\in\Z} \left( \ee^{-4\pi^2 (k+a-\alpha)^2 t}\left(
\frac{2(k+a)-\alpha}{4\pi^2\alpha t}-\frac{1}{(4\pi^2\alpha t)^2}\right) +  \frac{\ee^{-4\pi^2 (k+a)^2 t}}{(4\pi^2\alpha t)^2}\right)\ .
\end{eqnarray*}
Since $E^\alpha  E^{-\alpha}=-E^{-\alpha}  E^\alpha$ we obtain
by some resummation
\begin{eqnarray*}
\lefteqn{-t^2 E^\alpha  E^{-\alpha} \Tr_s \int_{\Delta^2}\ee^{t\sigma_0\partial_x^2} R_\alpha \ee^{t\sigma_1\partial_x^2} R_{-\alpha} \ee^{t\sigma_2\partial_x^2} d\sigma}\\
&=& \frac{ 2t}{\alpha}
E^\alpha E^{-\alpha}  \sum_{k\in\Z}(k+a) \ee^{-4\pi^2 (k+a)^2 t}\ .
\end{eqnarray*}
Using 
$$ t^{-2}(1+2t\frac{d}{dt})=(5+2t\frac{d}{dt})t^{-2}$$
we obtain
\begin{eqnarray*}
  \lefteqn{- E^\alpha  E^{-\alpha}t^{-1}\Psi_t(1+2t\frac{d}{dt}) t^2 \Tr_s \int_{\Delta^2}\ee^{t\sigma_0\partial_x^2} R_\alpha \ee^{t\sigma_1\partial_x^2} R_{-\alpha} \ee^{t\sigma_2\partial_x^2} d\sigma}\\&&
  - E^\alpha  E^{-\alpha}t^{-2}(1+2t\frac{d}{dt}) t^2 \Tr_s\int_{\Delta^2}\ee^{t\sigma_0\partial_x^2} R_\alpha \ee^{t\sigma_1\partial_x^2} R_{-\alpha} \ee^{t\sigma_2\partial_x^2} d\sigma\\
&=& (5+2t\frac{d}{dt})\frac{2}{t\alpha}
E^\alpha E^{-\alpha}\sum_{k\in\Z}(k+a) \ee^{-4\pi^2 (k+a)^2 t}\ .
\end{eqnarray*}
We employ
$$ \ee^{-4\pi^2 (k+a)^2 t}=\frac{1}{2\pi^{1/2}}\int_{-\infty}^\infty
\ee^{-z^2/4}\ee^{-2\pi\imath (k+a) t^{1/2} z} dz$$
in order to write
\begin{eqnarray*}
\lefteqn{(5+2t\frac{d}{dt})\frac{2}{t\alpha}\sum_{k\in\Z}(k+a) \ee^{-4\pi^2 (k+a)^2 t}}\\&=&(5+2t\frac{d}{dt})\frac{2}{t\alpha}\sum_{k\in\Z}
 \frac{1}{2\pi^{1/2}}\int_{-\infty}^\infty
\ee^{-z^2/4}\frac{-1}{2\pi\imath t^{1/2}}\frac{d}{dz}\ee^{-2\pi\imath (k+a) t^{1/2} z} dz\\
&=&(5+2t\frac{d}{dt}) \frac{\imath}{4 t^{3/2}\alpha\pi^{3/2}} 
\sum_{k\in\Z}\int_{-\infty}^\infty z \ee^{-z^2/4}  \ee^{-2\pi\imath (k+a) t^{1/2} z} dz\\
&=&(5+2t\frac{d}{dt}) \frac{\imath}{4 t^{3/2}\alpha\pi^{3/2}}\int_{-\infty}^\infty
\ee^{-2\pi\imath a t^{1/2} z}\sum_{k\in\Z}\ee^{-2\pi\imath t^{1/2} z k}
 z \ee^{-z^2/4}  dz\\
&=&(5+2t\frac{d}{dt}) \frac{\imath}{4 t^{3/2}\alpha\pi^{3/2}}\int_{-\infty}^\infty
\ee^{-2\pi\imath a t^{1/2} z}\sum_{m\in\Z}\delta(t^{1/2} z-m)
 z \ee^{-z^2/4}  dz\\
&=&(5+2t\frac{d}{dt})  \frac{\imath}{4 t^{5/2}\alpha\pi^{3/2}} \sum_{m\in\Z} m \ee^{-m^2/4t} \ee^{-2\pi\imath a m}\\
&=& \frac{\imath}{2t^{3/2}\alpha\pi^{3/2}} \sum_{m\in\Z} m \frac{d}{dt} \ee^{-m^2/4t} \ee^{-2\pi\imath a m}\\
&=& \frac{\imath}{8 t^{7/2}\alpha\pi^{3/2}} \sum_{m\in\Z} m^3 \ee^{-m^2/4t} \ee^{-2\pi\imath a m}\ .
\end{eqnarray*}
Integrating from $0$ to $\infty$ with respect to $t$ and substituting
$t=m^2/z$ we obtain
\begin{eqnarray*}
\lefteqn{\int_0^\infty \frac{\imath}{8 t^{7/2}\alpha\pi^{3/2}} \sum_{m\in\Z} m^3 \ee^{-m^2/4t} \ee^{-2\pi\imath a m}dt}\\&=&\frac{\imath}{8\alpha\pi^{3/2}} \sum_{m\in\Z\setminus \{0\}} \sign(m)\frac{1}{m^2} \ee^{-2\pi\imath a m} \int_0^\infty z^{3/2}\ee^{-z^2/4} dz\\
&=&\frac{1}{4\alpha\pi^{3/2}}  \int_0^\infty
z^{3/2}\ee^{-z^2/4} dz \sum_{m\ge 1} \frac{1}{m^2} \sin(2\pi a m)\ .
\end{eqnarray*}
Note that
\begin{eqnarray*}
 \int_0^\infty
z^{3/2}\ee^{-z^2/4} dz&=&
 \frac{1}{2}\int_0^\infty
z^{1/2}\ee^{-z^2/4} 2zdz\\
&=& \frac{1}{2} \int_0^\infty u^{1/4}\ee^{-u/4} du\\
&=&2^{3/2} \int_0^\infty v^{1/4}\ee^{-v} dv\\
&=&2^{3/2} \int_0^\infty v^{5/4-1}\ee^{-v} dv\\
&=&2^{3/2} \Gamma(5/4)\ .\end{eqnarray*}
We conclude 
\begin{lem}\label{uli1}
Let $V$ be spanned by $h_\alpha$, $|\alpha|\le R$. Then
$$ (T_2)_{|V}=\frac{1}{2^{1/2} \pi^{5/2}\imath} \Gamma(5/4)\left( \sum_{m\ge 1} \frac{1}{m^2} \sin(2\pi a m)\right)\sum_{\alpha=1}^R \alpha^{-1}E^\alpha E^{-\alpha}\ .$$
\end{lem}

Next we compute $\cD_* \cT_2(S^1,\cF)$.
For $k\in\Z$ let $X_k(x):=\ee^{2\pi\imath k x} \partial_x\in \cX(S^1)$. Then $h_{X_k}(x)= -4\pi\imath k \ee^{2\pi\imath x k} dx^2$.
By Lemma \ref{r33} and \ref{uli1}  we have 
\begin{eqnarray*}
\cD c_{T_2}(X_k,X_h)&=& T_2(h_{X_k},h_{X_h})\\
 &=& \frac{1}{2^{1/2}\pi^{5/2}\imath k} \Gamma(5/4)(-4\pi\imath k \ee^{2\pi\imath x k}) (-4\pi\imath h \ee^{2\pi\imath x h}) \delta_{k+h} \sum_{m\ge 1} \frac{1}{m^2} \sin(2\pi a m)\\
&=& 2^{7/2}\pi^{-1/2}\imath  k\Gamma(5/4)   \delta_{k+h}\sum_{m\ge 1} \frac{1}{m^2} \sin(2\pi a m)
\end{eqnarray*}
Define $u\in C^1(\cX(S^1))$ by
$u(f \partial_x):=\int_{S^1} f(x) dx$.
Then $du(f_1 \partial_x,f_2 \partial_x)=u([f_1 \partial_x,f_2 \partial_x])$. In particular,
\begin{eqnarray*}
d u(X_k,X_h)&=&2\pi\imath (h-k) \int_{S^1}  \ee^{2\pi\imath (k+h)} dx\\
&=&2\pi\imath (h-k) \delta_{h+k}\\
&=&-4\pi\imath k \delta_{h+k}\ .
\end{eqnarray*}
We conclude that
$$\cD c_{T_2}= - 2^{3/2}\pi^{3/2} \Gamma(5/4) \left( \sum_{m\ge 1} \frac{1}{m^2} \sin(2\pi a m)\right)  d u\ .$$
Thus we have shown the following
\begin{lem}
$$\cD_*\cT_2(S^1,\cF)=0\ .$$
\end{lem}
 
Finally we compute $\cT_2(S^1,\cF)\in H^2_c(\Diff(S^1)^0,\R)$.
By \cite{fuks84}, Thm. 3.4.4, the ring $H^*_c(\Diff(S^1)^0,\R)$
is generated by two classes $\alpha,\beta\in H^2_c(\Diff(S^1)^0,\R)$,
where $\beta^2=0$, $\cD_*\beta=0$, and $\cD_*\alpha\not=0$.
Thus $ \cT_2(S^1,\cF) = c \beta$ for some constant $c\in\R$.
It remains to describe $\beta$ and to determine the constant $c$.
As explained in \cite{fuks84},  there is a commutative diagram
$$\begin{array}{cccccccc}
 0   & \rightarrow  &  H^1_{top}(\Diff(S^1)^0,\R)  &  \stackrel{b}{\rightarrow}  &  H^2_c(\Diff(S^1)^0,\R)  &  \rightarrow  &  H^2(\cX(S^1),\R)  &  \rightarrow \\
     &              &    a\:\downarrow \cong            &               &
w\:\:\downarrow            &               &                    &\\
  0   &  \rightarrow  &  H^1_{top}(SL(2,\R),\R)  &  \stackrel{v}{\rightarrow}  &  H^2_c( SL(2,\R),\R)  &  \rightarrow  & 0  &  \rightarrow \end{array} \ ,$$
where the vertical maps are induced by the inclusion $SL(2,\R)\hookrightarrow \Diff(S^1)^0$ which is a homotopy equivalence of topological spaces,
$H^*_{top}$ denotes cohomology of topological spaces, 
and we have used that $H^1(sl(2,\R),\R)=0$, $H^2(sl(2,\R),\R)=0$, $H^1(\cX(S^1),\R)=0$.
 
The Cartan decomposition gives a decomposition $SL(2,\R)=SO(2)\times \R^2$ as topological space, where  $SO(2)\cong S^1$ is  a maximal compact subgroup.
Let $f:SL(2,\R)\rightarrow SO(2)$ be projection onto the first factor.
 
We fix the orientation of $SO(2)$ such that $$\left(\begin{array}{cc} 0&-1\\1&0\end{array}\right)\in so(2)$$ points into the positive direction.
Let $\delta\in  H^1_{top}(SL(2,\R),\R)$ be such that $\langle \delta,[SO(2)\times\{0\}]\rangle=1$.
Then $\beta:=b\circ a^{-1} (\delta) $.

We determine the constant $c$ using the identity $w(\cT_2(S^1,\cF))= c v(\delta)$.
By \cite{fuks84}, Thm. 3.4.3., there is a unique sheet of 
$\Arg$  such that
$c^\prime:SL(2,\R)\times SL(2,\R)\rightarrow \R$ given
by $c^\prime(g,h):=\frac{1}{2\pi}\Arg(f(gh)-f(h)-f(g))$
satisfies $c^\prime(1,1)=0$. The function 
$c^\prime$ represents $v(\delta)$ using the description
\cite{fuks84}, Ch. 4.1 (II), of $H^2_c(SL(2,\R),\R)$.
As explained in \cite{fuks84}, p. 264, the class $v(\delta)$
is also represented by the cocylce 
$(g,h)\mapsto \frac{1}{2\pi}\vol_{H^2}(o, go,gho)$, where $o$ is the origin
in the hyperbolic plane $H^2=SL(2,\R)/SO(2)$ of constant sectional curvature $-1$ and $\vol_{H^2}(o, go,gho)$
denotes the oriented volume of the geodesic triangle $\Delta(g,h)$
spanned by $o, go,gho\in H^2$. 

Note that
$w(\cT_2(S^1,\cF))$ is represented by $c^\prime_{T_2}(g,h):=c_{T_2}(1,g,gh)$.
Since $SO(2)$ stabilizes the metric $dx^2$,
the inclusion $SL(2,\R)\hookrightarrow \Diff(S^1)^0$
induces an inclusion $i:H^2\hookrightarrow \Met_1(S^1)$
such that $i(o)=dx^2$.

Let $I$ denote the unit interval $[0,1]$. 
If $g\in SL(2,\R)$, then we define
$j(g):I^2\mapsto \Met(S^1)$ by
$j(g)(s,t)=s i(\gamma_t) + (1-s)(tdx^2+(1-t) w(g) dx^2)$,
where $\gamma:I\mapsto H^2$ is the geodesic path joining 
$o$ and $go$. Furthermore we define the cochain $u\in C^1_c(\Diff(S^1)^0)$ by $u(g):=\int_{I^2} j(g)^*T_2$.
Then we have $c^\prime_{T_2}(g,h)-\int_{\Delta(g,h)} i^* T_2 = 
u(gh)-u(g)-u(h)$. It follows that
$w(\cT_2(S^1,\cF))$ is represented by the cocycle $(g,h)\mapsto \int_{\Delta(g,h)} i^* T_2$. Since $i$ is $SL(2,\R)$-equivariant and $T_2$ is $\Diff(S^1)^0$-invariant, $i^* T_2$ is $SL(2,\R)$-invariant
and hence proportional to the volume form $\vol_{H^2}$,
thus $i^* T_2 = \frac{c}{2\pi} \vol_{H^2}$.

Let
$$A:=\left(\begin{array}{cc} 1&0\\ 0&-1\end{array}\right),\quad N:=\left(\begin{array}{cc} 0&1\\ 0&0\end{array}\right)\ .$$
Let $A^\sharp, N^\sharp \in \cX(H^2)$ denote the corresponding
fundamental vector fields.  
Then $$ \vol_{H^2}(A^\sharp(o),N^\sharp(o))=-2\ .$$

Let $A^*, N^*\in \cX(S^1)$ denote the fundamental vector fields
corresponding to $A,N$. Then we have
$$A^*(x)=\frac{1}{2\pi\imath}(\ee^{2\pi\imath x}-\ee^{-2\pi\imath x})\partial_x,\quad
N^*(x)=\frac{1}{4\pi}(\ee^{2\pi\imath x}+\ee^{-2\pi\imath x}+2)\partial_x\ .$$
We have
\begin{eqnarray*}
i^*T_2(A^\sharp,N^\sharp)&=& 
\frac{d}{ds}_{|s=0} \frac{d}{dt}_{|t=0}c_{T_2}(1,\ee^{tA},\ee^{tA}\ee^{sN})\\
&=&\cD c_{T_2}(A^*,N^*)\\
&=&-\frac{1}{2^{5/2}\pi^{9/2}}     \Gamma(5/4) \sum_{m\ge 1} \frac{1}{m^2} \sin(2\pi a m)\ .
\end{eqnarray*}
We conclude that
$$w(\cT_2(S^1,\cF))=-\frac{1}{2^{5/2}\pi^{7/2}}     \Gamma(5/4) \sum_{m\ge 1} \frac{1}{m^2} \sin(2\pi a m)v(\delta)\ .$$
We have shown the following 
\begin{lem}
$$\cT_2(S^1,\cF)=-\frac{1}{2^{5/2}\pi^{7/2}}     \Gamma(5/4) \left( \sum_{m\ge 1} \frac{1}{m^2} \sin(2\pi a m)\right) \beta\ .$$
\end{lem}
 
\section{Concluding remarks}

\begin{enumerate}
\item The one-dimensional example shows that $\cT(M,\cF)$
is nontrivial in general.
\item In \cite{bott77} Bott gave a construction of cocyles
for $\Diff(M)$ given by integration of locally
computable quantities. In our example the class $\alpha$ 
can be represented in this way. Since $\cT(S^1,\cF)$
depends non-trivially on $\cF$ there is no local representation
for $\cT(M,\cF)$ in general.
\item Is there any easy way to compute $\cT_{2j}(S^1,\cF)$
for $j\ge 2$?
\item Can $\cD_* \cT(M,\cF)$ be non-trivial?
\item Give a differential topological interpretation of $\cT(M,\cF)$?
\end{enumerate}

\bibliographystyle{plain}

\end{document}